# Networks and their degree distribution, leading to a new concept of small worlds

Leo Egghe

Hasselt University, Hasselt, Belgium

E-mail: leo.egghe@uhasselt.be

ORCID: 0000-0001-8419-2932

## Abstract

The degree distribution, referred to as the delta-sequence of a network is studied. Using the non-normalized Lorenz curve, we apply a generalized form of the classical majorization partial order.

Next, we introduce a new class of small worlds, namely those based on degree centralities of networks. Similar to a previous study, small worlds are defined as sequences of networks with certain limiting properties. We distinguish between three types of small worlds: those based on the highest degree, those based on the average degree, and those based on the median degree. We show that these new classes of small worlds are different from those introduced previously based on the diameter of the network or the average and median distance between nodes. However, there exist sequences of networks that qualify as small worlds in both senses of the word, with stars being an example. Our approach enables the comparison of two networks with an equal number of nodes in terms of their "small-worldliness".

Finally, we introduced neighboring arrays based on the degrees of the zeroth and first-order neighbors and proved that for trees, equal neighboring arrays lead to equal delta-arrays.

## 1. Introduction

Consider an undirected network or graph $G = (V,E)$, where $V$ denotes the set of vertices or nodes, and $E$ denotes the set of edges or links. In this text, the terms graph and network refer to the same mathematical concept and are used interchangeably. A path of length n is a sequence of vertices $(v_0, ...v_k, v_{k+1}, ..., v_n)$ such that $\{v_0, ..., v_{n-1}\}$ and $\{v_1, ..., v_n\}$ are sets (being sets each consist of different elements) and for $k = 0,..., n-1$, $v_k$ is adjacent to $v_{k+1}$. A cycle is a path for which the starting point $v_0$ coincides with the endpoint $v_n$. A graph is connected if there exists (at least one) path between any two vertices. If $\#V = N$, then the degree centrality of node i, i = 1, ..., N, i.e., the number of edges connected to node i, is denoted as $\delta_i$. In this article we always assume that G is connected, hence all degree centralities are strictly larger than zero. As there is no natural order among the nodes in a network we assume that these values are ranked in decreasing order.

Notation

The array of degree centralities of a network G with N nodes is denoted as

$$\Delta_G = (\delta_1(G), \delta_2(G), \dots, \delta_N(G)). \qquad (1)$$

We will informally refer to such an array as a delta-sequence, consisting of delta-values. Clearly, $\sum_{i=1}^{N} \delta_i = 2\,(\#E)$, a notion which is known as the total degree of the network. It is easy to see that $2(N-1) \leq \sum_{i=1}^{N} \delta_i \leq N(N-1)$. The lower bound is obtained e.g., for a chain consisting of N nodes, see further, while the upper bound is obtained for a complete graph where each node is connected to each other node.

As real-world networks are often dynamic we will work, as we did in a previous study on small-world networks (Egghe & Rousseau, 2024), within the context of a sequence of finite node sets $(V_N)_{N \in \mathbb{N}}$ of networks $(G_N)_{N \in \mathbb{N}}$, i.e., $\forall N \in \mathbb{N}$: $\#V_N = N$ Of course, also all edge sets $E_N$ are finite.

Before moving on to examples and theory we recall the following definitions.

## 1.1 Definition:  Free or unrooted tree (Knuth, 1973, p. 363)

A free or unrooted tree is a connected graph with no cycles. Equivalently it is a connected graph such that removing any edge makes it disconnected. Another equivalent definition states that if v and v′ are different vertices, then there exists exactly one path from v to v′. As we will never use the notion of the root of a tree, we will just use the term 'tree' for "free tree".

## 1.2 Definition. Isomorphic graphs

Two graphs G and G′ are isomorphic if there exists a bijection f between the vertices of G and G′ such that there is an edge between vertices u and v in G if and only if there is an edge between the vertices f(u) and f(v) in G′.

## 1.3 Definition: Spanning tree of a connected graph

A spanning tree of an N-node connected graph is a set of N-1 edges that connects all nodes of the network and contains no cycles. A graph may have different (non-isomorphic) spanning trees.

## 2. Examples of networks and their delta-sequences

### 2.1 The complete network on N nodes

The delta-sequence of an N-node complete network is

$$\Delta_G = (\underbrace{N-1, \ldots, N-1}_{N\ times})$$

In this case $\sum_{i=1}^{N} \delta_i$ = N(N-1) and its diameter is 1.

### 2.2 The N-star

The N-star consists of a central node and N-1 peripheral nodes, each with one link, namely to the center. Then

$$\Delta_G = \left(N-1, \ \underbrace{1, \ldots, 1}_{N-1-times}\right)$$

Its sum is 2(N-1) and its diameter is 2.

### 2.3 The N-polygon

The N-polygon consists of N nodes forming a simple path of different nodes, that links to its starting point. Then

$$\Delta_G = \left( \underbrace{2, \dots, 2}_{N \; times} \right)$$

Its sum is 2N and its diameter is N/2 for even N and (N-1)/2 for odd N.

2.4 The N-chain

The N-chain consists of one path of N different nodes. Then

$$\Delta_G = \left( \underbrace{2, \dots, 2}_{N-2 \; times}, 1,1 \right)$$

Its sum is 2(N-1) and its diameter is N-1.

2.5 Trees

It is obvious that there is no delta-sequence applicable to all trees, but we do have the following lemma, to which we will refer as Knuth's lemma.

Lemma (Knuth, 1976, p. 363)

An N-node connected network is a tree if and only if it has N-1 edges and hence a total degree equal to 2(N-1).

We note that stars and chains are all special trees.

## 3. The delta-sequence and a generalized majorization partial order

3.1 The standard Lorenz curve and Gini index

As the delta-sequence does not have a fixed sum we first consider its (standard) Lorenz curve. The highest Lorenz curve for a network with N nodes is obtained for the star. It always starts by connecting the origin to the point with coordinates $\left(\frac{1}{N}, \frac{1}{2}\right)$. In general, its standard Gini index (Rousseau et al., 2018, formula (4.19)) is (N-2)/2N with a limiting value of 0.5. The lowest Lorenz curve, with Gini index zero, is obtained for a delta-sequence consisting of the same numbers, such as for any complete network, but also for any polygon.

## 3.2 The non-normalized Lorenz curve

In (Egghe & Rousseau, 2023a) we used the so-called non-normalized Lorenz curve in a continuous context. This study and its follow-up (Egghe & Rousseau, 2023b) led to a rigorous definition of the notion of global impact. We would say, based on (Egghe & Rousseau, 2023b), that the notion of majorization in a network does not only depend on the number of links, but also on their concentration.

In the discrete context, the non-normalized Lorenz curve is defined as follows.

Definition: Non-normalized Lorenz curves

Let $X = (x_1, x_2, \dots, x_N)$ be a decreasing N-array of non-negative real numbers, then the corresponding non-normalized Lorenz curve is the polygonal line connecting the origin (0,0) with the points $\left(j, \sum_{k=1}^{j} x_j\right)$, j= 1, …, N. This curve ends at the point with coordinates $\left(N, \sum_{k=1}^{N} x_j\right)$.

Definition: The non-normalized (or generalized) majorization order for N-arrays

If X and Y are decreasing N-arrays of non-negative real numbers, then X is majorized by Y, denoted as X -< Y if

$$\forall j, j = 1, \dots, N: \sum_{k=1}^{j} x_j \leq \sum_{k=1}^{j} y_j \tag{2}$$

The -< relation is only a partial order as non-normalized Lorenz curves (just like standard Lorenz curves) may intersect, see further. If there exists at least one j such that the inequality is strict, we say that X is strictly majorized by Y. If $\forall j, j = 1, \dots, N; x_j \leq y_j$ then obviously X -< Y, but the opposite relation does not hold.

Definition. Acceptable measures

If **X** denotes the set of all decreasing N-arrays of non-negative real numbers, then a function m: **X** $\rightarrow \mathbb{R}^+$ is an acceptable measure for the relation -< if X -< Y implies that m(X) ≤ m(Y).

It is important to note that the definitions of non-normalized Lorenz curves and in particular the notion of the generalized Lorenz majorization order, denoted as -<, and the corresponding acceptable measures are generally applicable to all decreasing N-array of non-negative real numbers.

Hence, the above definitions can be applied to the set of delta-sequences and the corresponding networks, leading to expressions such as $\Delta_G$ -< $\Delta_H$ for two N-node networks, but also to the gamma-sequences, introduced later in this text. We already note that if T(G) denotes a spanning tree of the network G, then T(G) -< G.

As an illustration of the importance of the generalized majorization order, we recall the following definition.

Definition: Network density (Wasserman & Faust, 1994)

The density D of an undirected network G with N nodes is defined as

$$D(G) = \frac{2.(\#E)}{N(N-1)} = \frac{\sum_{i=1}^{N} \delta_i}{N(N-1)} \tag{3}$$

Clearly, network density is just a normalized total degree. Two N-networks with the same density D have non-normalized Lorenz curves with the same endpoint, but D does not say anything about the exact relation between the two non-normalized Lorenz curves. In this sense, the majorization partial order applied to delta-sequences refines the notion of network density.

Let now $\Delta_G$ and $\Delta_H$ be the degree sequences of the N-node networks G and H, then the following theorem holds.

3.4 Theorem 1

(i) $\Delta_G$ -< $\Delta_H \Rightarrow \delta_1(G) \leq \delta_1(H)$

(ii) $\Delta_G$ -< $\Delta_H \Rightarrow \frac{1}{N} \sum_{j=1}^{N} \delta_j(G) \leq \frac{1}{N} \sum_{j=1}^{N} \delta_j(H)$

(iii) $\Delta_G$ -< $\Delta_H \nRightarrow$ Md($\Delta_G$) ≤ Md($\Delta_H$) and neither does it imply that Md($\Delta_G$) ≥ Md($\Delta_H$), where Md stands for the median of a sequence.

(iv) The reverse implications do not hold

Proof. (i) and (ii) follow trivially from the definition of -<.

(iii) We provide two counterexamples (N=5)

The 5-chain G has a degree sequence $\Delta_G$ = (2,2,2,1,1) and H (see Fig. 1) has a degree sequence $\Delta_H$ = (3,2,1,1,1). Then $\Delta_G$ -< $\Delta_H$ but M$d$($\Delta_G$) = 2 > Md($\Delta_H$) = 1.

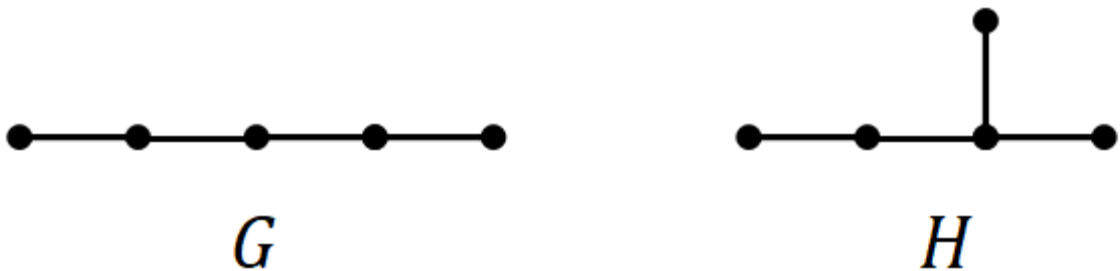

Fig. 1. Networks G and H illustrating part (iii)

For $G_1$ and $H_1$ (see Fig. 2) we have $\Delta_{G_1}$ = (3,3,2,2,2) and $\Delta_{H_1}$ = (4,3,3,2,2). Then $\Delta_{G_1}$ -< $\Delta_{H_1}$ and M$d$($\Delta_{G_1}$) = 2 < M$d$($\Delta_{H_1}$) = 3.

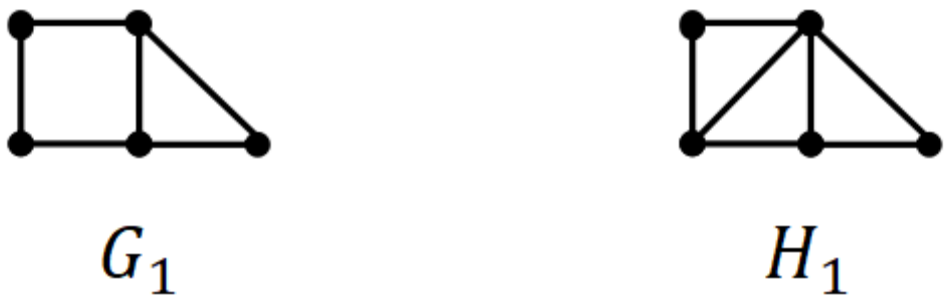

Fig. 2 Networks $G_1$ and $H_1$ illustrating part (iii)

(iv) For the opposite of case (i) we consider the networks $H_1$ and $H_2$ (see Fig. 3) with $\Delta_{H_1} = (5,3,3,3,3,3)$ and $\Delta_{H_2} = (4,4,4,3,3,2)$. Then neither $\Delta_{H_1}$ -< $\Delta_{H_2}$ , nor $\Delta_{H_2}$ -< $\Delta_{H_1}$ , illustrating that -< is a partial, not a complete, order. Moreover, $\delta_1(H_2) = 4 < \delta_1(H_1) = 5$ , $\overline{\delta(H_1)} = \overline{\delta(H_2)} = \frac{20}{6}$ and $Md(\Delta_{H_2}) = \frac{7}{2} > Md(\Delta_{H_1}) = 3$ □

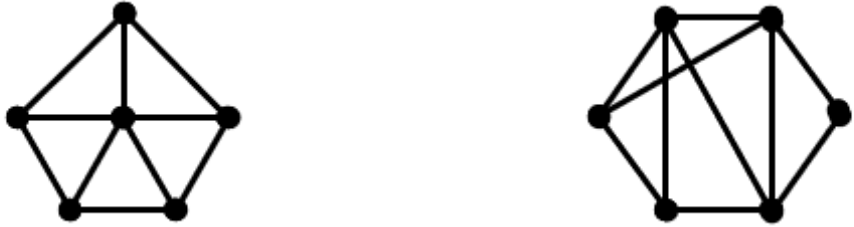

Fig.3 The networks $H_1$ (left) and $H_2$ (right) used in part (iv)

From our earlier investigation (Egghe & Rousseau, 2023a), we know that the following are acceptable measures for generalized majorization among delta-sequences:

A. The Gini index: $\text{Gini}(\Delta_G) = \sum_{i=1}^{N}\left(\sum_{j=1}^{i} \delta_j\right)$

B. The entropy or Theil measure: $\text{Th}(\Delta_G) = \sum_{j=1}^{N} \delta_j \ln(\delta_j)$

C. The power measure: $P(\Delta_G) = \sum_{j=1}^{N} \delta_j^p, p > 1$

## 3.5 The network with the lowest non-normalized Lorenz curve.

Theorem 2. An N-node chain is the lowest connected network in the generalized majorization partial order.

Proof. By Knuth's Lemma, the endpoint of the generalized Lorenz curve of any network that is not a tree is situated strictly above that of a tree. Hence, the lowest possible network must be a tree. Among all trees, the N-chain has the lowest generalized Lorenz curve.

Remark. We note that it is not even possible for a general network to have a generalized Lorenz curve that at any place is situated below that of an N-chain. Indeed, when the endpoint is fixed, then the lowest generalized Lorenz curve is the one whose classical Lorenz curve is the diagonal. As the lowest possible endpoint is 2N, this corresponds e.g., to the N-polygon, with delta-array (2,2,..., 2). Its cumulative array is (2,4,6,...., 2N-4, 2N-2, 2N). Yet, the corresponding array for the N-chain is (2,4,6,...., 2N-2, 2N-1,2N), showing that it is not possible to be situated (locally) strictly under the generalized Lorenz curve of the N-chain.

Remark further that the largest generalized Lorenz curve of an N-node network is the one corresponding to the N-complete network.

3.6 Examples of non-comparable networks

We already remarked that $\prec$ is only a partial order, implying that some N-node networks are not comparable. Theorem 1

already provided some examples. Here we provide some more examples.

a) Case N = 5

Consider the networks in Fig. 4

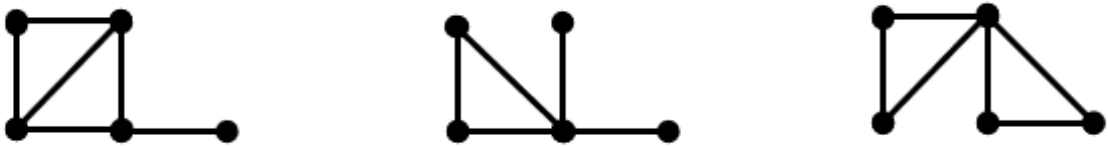

Fig.4. Incomparable 5-node networks: $G_1$, $G_2$ and $G_3$

For the first and the second network, we have delta-sequences (3,3,3,2,1) and (4,2,2,1,1). Hence $\Delta_{G_1} \nprec \Delta_{G_2}$ and $\Delta_{G_2} \nprec \Delta_{G_1}$. These networks have a different number of links. Consider now $G_1$ and $G_3$ with the same number of links. The delta-sequence of $G_3$ is (4,2,2,2,2). Then, clearly $\Delta_{G_1} \nprec \Delta_{G_3}$ and $\Delta_{G_3} \nprec \Delta_{G_1}$. We further note that $D(G_1) > D(G_2)$, $\delta_1(G_1) < \delta_1(G_2)$ and $\frac{1}{N} \sum_{j=1}^{N} \delta_j(G_1) > \frac{1}{N} \sum_{j=1}^{N} \delta_j(G_2)$ but $\Delta_{G_1} \nprec \Delta_{G_2}$ and $\Delta_{G_2} \nprec \Delta_{G_1}$.

b) Case N = 6 and beyond.

Consider again the networks shown in Fig.3.

The delta-sequences of $H_1$ and $H_2$ are respectively (5,3,3,3,3,3) and (4,4,4,3,3,2). Their sums are equal to 20. Yet, $\Delta_{H_1} \nprec \Delta_{H_2}$ and $\Delta_{H_2} \nprec \Delta_{H_1}$. If we add chains of equal length to a node with degree 3, we may obtain incomparable networks with any larger number of nodes.

c) The cases 2 ≤ N ≤ 4.

i) There is only one network with N = 2, hence any two networks are comparable.

ii) The case N=3. Then there are only 2 non-isomorphic networks, shown in Fig.5.

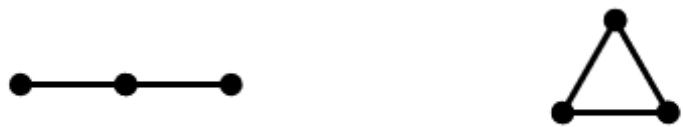

Fig. 5. The two non-isomorphic networks with three nodes.

The delta-sequence of the chain on the left is (2,1,1), while the delta-sequence of the polygon on the right is (2,2,2). Clearly, the chain is strictly smaller than the polygon.

iii) The case N = 4. There are six non-isomorphic connected 4-node networks. See Fig. 6.

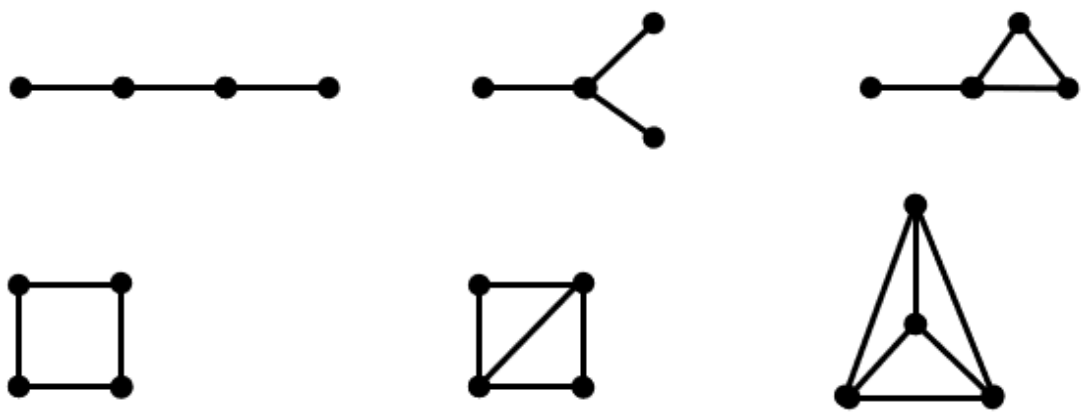

Fig. 6. The six non-isomorphic networks with degree 4.

Referring to these networks, from left to right and from the first row to the second, as $G_1$, $G_2$, $G_3$, $G_4$, $G_5$, and $G_6$ we obtain the following delta-sequences:

$\Delta_{G_1} = (2,2,1,1)$; $\Delta_{G_2} = (3,1,1,1)$; $\Delta_{G_3} = (3,2,2,1)$; $\Delta_{G_4} = (2,2,2,2)$; $\Delta_{G_5} = (3,3,2,2)$; and $\Delta_{G_6} = (3,3,3,3)$. The corresponding cumulative distributions are: $(2,4,5,6)$, $(3,4,5,6)$, $(3,5,7,9)$, $(2,4,6,8)$, $(3,6,8,10)$ and $(3,6,9,12)$. Hence we have the following relations between these networks (Fig. 7).

$$\Delta_{G_1} -< \begin{Bmatrix} \Delta_{G_2} \\ \Delta_{G_4} \end{Bmatrix} -< \Delta_{G_3} -< \Delta_{G_5} -< \Delta_{G_6}$$

Fig.7 Relations between delta-sequences of networks with 4 nodes

We see that only $G_2$ and $G_4$ are incomparable. We note that the (generalized) Gini-indices for these networks are: Gini($G_1$) = 17, Gini($G_2$) = 18, Gini($G_3$) = 24, Gini($G_4$) = 20, Gini($G_5$) = 27 and Gini($G_6$) = 30, illustrating the fact that the Gini index is an acceptable measure for the generalized majorization order.

## 4. Degree sequences and small worlds

### 4.1 Introduction to this section

In this section, we introduce a new class of small worlds, namely those based on degree centralities of networks. Similar to a previous study, small worlds are defined as sequences of networks with certain limiting properties. We distinguish between different types of small worlds and show that these new classes of small worlds are different from those introduced previously (Egghe & Rousseau, 2024). However, there exist sequences of networks that qualify as small worlds in both senses of the word, with stars being an example. Our approach enables the comparison of two networks with an equal number of nodes in terms of their "small-worldliness".

## 4.2. Definitions of small worlds derived from the degree distribution

Because we will define here small worlds derived from the degree distribution we will use the abbreviation DSW.

### 4.2.1 Small worlds derived from the largest degree (DSWL)

If $\Delta_{G_N} = (\delta_1(G_N), \delta_2(G_N), \ldots, \delta_N(G_N))$, $N \in \mathbb{N}$, is the delta-sequence of $G_N$, then $(G_N)_{N \in \mathbb{N}}$ is a degree small world based on the largest degree if

$$\lim_{N \to +\infty} \frac{\delta_1(G_N)}{\ln(N)} = +\infty \qquad (4)$$

We informally say that $(G_N)_{N \in \mathbb{N}}$ is DSWL.

### 4.2.2 Small worlds derived from the average degree (DSWA)

If $\overline{\delta(G_N)}, N \in \mathbb{N}$, denotes the average degree in network $G_N$,

$$\overline{\delta(G_N)} = \frac{1}{N} \sum_{j=1}^{N} \delta_j(G_N) \qquad (5)$$

then $(G_N)_{N \in \mathbb{N}}$ is a degree small world based on the average degree if

$$\lim_{N \to +\infty} \frac{\overline{\delta(G_N)}}{\ln(N)} = +\infty \qquad (6)$$

In the same vein as above we say that $(G_N)_{N \in \mathbb{N}}$ is DSWA.

### 4.2.3 Small worlds derived from the median degree (DSWMd)

If $Md(G_N), N \in \mathbb{N}$, denotes the median degree of the network $G_N$, then $(G_N)_{N \in \mathbb{N}}$ is a degree small world based on the median degree if

$$\lim_{N \to +\infty} \frac{Md(G_N)}{\ln(N)} = +\infty \qquad (7)$$

We say that $(G_N)_{N \in \mathbb{N}}$ is DSWMd.

### 4.3. Proposition 1

a) If $(G_N)_{N \in \mathbb{N}}$ is DSWMd, then $(G_N)_{N \in \mathbb{N}}$ is DSWA

b) If $(G_N)_{N \in \mathbb{N}}$ is DSWA, then $(G_N)_{N \in \mathbb{N}}$ is DSWL

c) If $(G_N)_{N \in \mathbb{N}}$ is DSWMd, then $(G_N)_{N \in \mathbb{N}}$ is DSWL

d) the reverse relations do not hold.

Proof. Implication a) follows from the Markov property (Chow & Teicher, 1978), which states that $Md(G_N) \leq 2\overline{\delta(G_N)}$. Hence, if $(G_N)_{N \in \mathbb{N}}$ is DSWMd, then $(G_N)_{N \in \mathbb{N}}$ is DSWA.

Implication b) follows from the fact that $\overline{\delta(G_N)} \leq \delta_1(G_N)$.

Implication c) follows immediately from implications a) and b), or by noticing that $Md(G_N) \leq \delta_1(G_N)$.

d) Consider the following sequence of star networks $(ST_N)_N$ (Fig. 8)

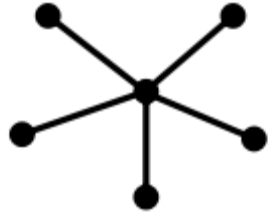

Fig. 8. Star network (illustrated for N=6)

Then $\Delta_{ST_N} = \left( N-1, \underbrace{1, ..., 1}_{(N-1) times} \right)$. Clearly $(ST_N)_N$ is DSWL, but not DSWA or DSWMd.

Finally, we have to show that DSWA does not imply DSWMd.

Consider Fig.9, to which we refer as an M-spider, denoted as $S_M$ (Fig.9 shows a spider with M =5). It consists of a complete M-node graph, where each node has an extra two links. Hence N= 3M.

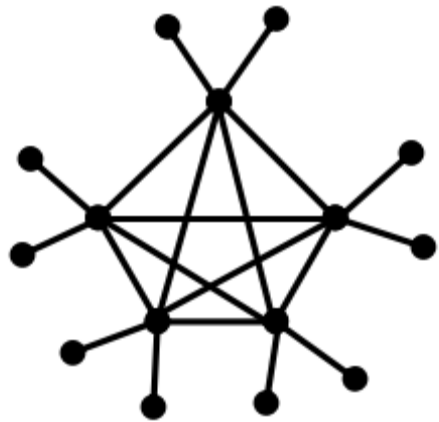

Fig.9. 5-spider

The delta-sequence of an M-spider is $\Delta_{S_M} = \left( \underbrace{M+1, \ldots, M+1}_{M \text{ times}}, \right.$

$\left. \underbrace{1, \ldots, 1}_{2M \text{ times}} \right)$. The average $\frac{1}{N} \sum_{j=1}^{N} \delta_j(S_M) = \frac{M^2+3M}{3M}$ and hence

$\lim_{N \to +\infty} \frac{\overline{\delta(S_M)}}{\ln(N)} = \lim_{M \to +\infty} \frac{\left(\frac{M}{3}\right)+1}{\ln(3M)} = +\infty$, showing that $(S_M)_M$ is DSWA.

As $Md(\Delta_{S_M}) = 1$, this shows that $(S_M)_M$ is not DSWMd.□

## 4.4 Examples

We have already considered the sequence of stars. Now, we have a look at complete graphs, chains, and polygons.

### 4.4.1. The sequence of complete graphs.

For each N, $\Delta_{G_N} = (\underbrace{N-1, \ldots, N-1}_{N \text{ times}})$. Hence complete graphs are DSWMd and hence also DSWA and DSWL.

### 4.4.2. The sequence of chains of length N, $(CH_N)_N$. Then

$\Delta_{CH_N} = \left( \underbrace{2, \ldots, 2}_{(N-2) \text{ times}}, 1, 1 \right)$. Chains are not small worlds based on their degree sequences.

### 4.4.3. Polygons $(Po_N)_N$. These are formed by connecting begin and end nodes of chains. Their delta-sequences are: $\Delta_{Po_N} = \left( \underbrace{2, \ldots, 2, 2}_{N \text{ times}} \right)$. We see that polygons too are not small worlds based on their degree sequences.

## 5. Another approach to small worlds

In a previous article (Egghe & Rousseau, 2024) we studied small worlds using the following definitions. We assume given a sequence $(\Omega_N)_{N \in \mathbb{N}}$ of finite sets and a distance function d defined on $\mathbf{\Omega} = \cup_{N \in \mathbb{N}} \Omega_N$. The term "small world" was used there for a sequence $(\Omega_N)_{N \in \mathbb{N}}$ of finite sets satisfying one of the properties defined below.

### 5.1 Small worlds based on the diameter (SWD)

If $d_N$, $N \in \mathbb{N}$, is the diameter of $\Omega_N$, defined as

$$d_N = max\{d(A,B); A,B \, \epsilon \, \Omega_N\} \qquad (8)$$

then $(\Omega_N)_{N \in \mathbb{N}}$ is a SWD if there exists a finite constant $C \geq 0$ such that

$$\lim_{N \to +\infty} \frac{d_N}{\ln(N)} = C \qquad (9)$$

Note that $d_N$ is short for $diam(\Omega_N)$.

### 5.2 Small worlds based on the average distance (SWA)

If $\mu_N$, $N \in \mathbb{N}$, denotes the average distance between two different elements in $\Omega_N$:

$$\mu_N = \frac{1}{N(N-1)} \sum_{\substack{A,B \, \epsilon \, \Omega_N \\ A \neq B}} d(A,B) \qquad (10)$$

then $(\Omega_N)_{N \in \mathbb{N}}$ is an SWA if there exists a finite number $C \geq 0$ such that

$$\lim_{N \to +\infty} \frac{\mu_N}{\ln(N)} = C \qquad (11)$$

### 5.3 Small worlds based on the median distance (SWMd)

If $Md_N, N \in \mathbb{N}$, denotes the median distance between two different elements in $\Omega_N$:

$$Md_N = median\{\{d(A,B); A, B \in \Omega_N, A \neq B\}\} \qquad (12)$$

then $(\Omega_N)_{N \in \mathbb{N}}$ is a SWMd if there exists a finite number $C \geq 0$ such that

$$\lim_{N \to +\infty} \frac{Md_N}{\ln(N)} = C \qquad (13)$$

Note that $\{\{...\}\}$ in (12) refers to a multiset (Rousseau et al., 5.13.1), i.e. a "set" in which elements may occur more than once. It is obvious that if a sequence of finite sets is an SWD $(\Omega_N)_{N \in \mathbb{N}}$ then it is also an SWA and an SWMd (Egghe & Rousseau, 2024, section 2.4).

## 6. The relation between SWs and DSWs

We recall from (Egghe, 2024) that if $\alpha_j$ (G), j= 1,..., N-1, denotes the number of times distance j (the shortest distance between two nodes) occurs in network G, then the array $AF_G = (\alpha_1(G), \alpha_2(G), ..., \alpha_{N-1}(G))$ is called the $\alpha -$ array of the network G.

In this section, we will find out if being an SW (Egghe & Rousseau, 2024) based on the so-called alpha-sequence, implies also being a DSW or vice versa. It will be shown that such implications do not exist, which implies that the notions of SW and DSW are different concepts.

Two sequences of relations are already known: SWD ⇒ SWA ⇒ SWMd, (Egghe & Rousseau, 2024) and DSWMd ⇒ DSWA ⇒

DSWL (see above). We will prove now that there is, in general, no relation between these sequences of implications.

## 6.1 Theorem 3

(i) SWD ⇏ DSWL

(ii) DSWMd ⇏ SWMd

Proof.(i). We first construct a network $\Omega_N$ for fixed N > 3. Consider a chain with $\lfloor \ln(N) \rfloor$ nodes (step 1). Each of these points has $\lfloor \ln(N) \rfloor$ descendants (see Fig. 10) (step 2).

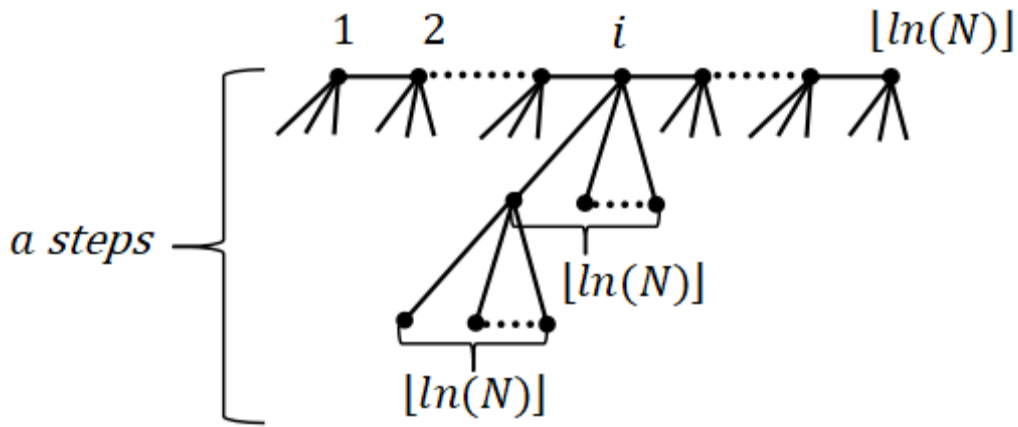

Fig.10 Sketch of the construction of a counterexample, used in Theorem 3 (i)

We continue this construction until at step a, we have $\sum_{k=1}^{a}(\lfloor \ln(N) \rfloor)^k \geq N$. We see that a is smaller than or equal to any b for which $(\lfloor \ln(N) \rfloor)^b \geq N$. Then the number b satisfies the inequality $b \geq \frac{\ln(N)}{\ln(\lfloor \ln(N) \rfloor)}$. We see that this network's diameter $d_N$ satisfies the inequality $d_N \leq \lfloor \ln(N) \rfloor + \frac{2\ln(N)}{\ln(\lfloor \ln(N) \rfloor)}$. Hence

$$\lim_{N \to \infty} \frac{d_N}{\ln(N)} < +\infty$$

which proves that $(\Omega_N)_{N \in \mathbb{N}}$ is SWD. However, $(\Omega_N)_{N \in \mathbb{N}}$ is not DSW because each delta-value is smaller than or equal to $\lfloor \ln(N) \rfloor + 2$, so that $\lim_{N \to \infty} \frac{\delta_1}{\ln(N)}$ cannot be equal to $+\infty$.

(ii). We will construct a kite consisting of an M-complete network and a tail consisting of M-1 nodes (Fig. 11). Hence N = 2M-1. For fixed N the delta-sequence of this kite, $K_N$, is

$$\Delta_{K_N} = \left( M, \underbrace{M-1, \dots, M-1}_{(M-1)\,times}, \underbrace{2, \dots, 2}_{(M-2)\,times}, 1 \right).$$ We see that $\sum_{j=1}^{N} \delta_j(K_n) =$

$M + (M-1)^2 + 2(M-2) + 1 = M^2 + M - 2$.

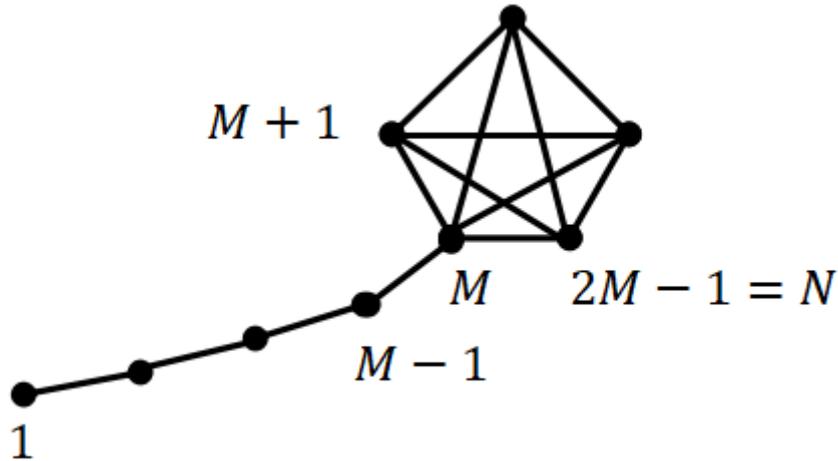

Fig. 11. Kite with N=2M-1 nodes

Now, the kite's alpha-sequence is: $AF(K_N) = \left( \left( \frac{M(M-1)}{2} \right) + M - \right.$

$1, 2M-3, 2M-4, \dots, M-1, \underbrace{0, \dots, 0}_{(N-M-1)\,times} \left. \right)$. One may check that

$\alpha_1(K_N) = \frac{M^2}{2} + \frac{M}{2} - 1 = \frac{1}{2}\sum_{j=1}^{N} \delta_j(K_N)$ and that $\sum_{j=1}^{N-1} \alpha_j(K_N) = \frac{N(N-1)}{2}$.

Now we see that $\forall i \in \{1, \ldots, N\}: \sum_{j=1}^{i} \alpha_j(K_N) = \frac{M(M-1)}{2} + (M-1) +$

$(2M-3) + (2M-4) + \cdots + 2M - (i+1) \quad = \quad \frac{M(M-1)}{2} + (M-1) +$

$2M(i-1) - (3 + 4 + \cdots + i + 1) \quad = \quad \frac{M(M-1)}{2} + (M-1) + 2M(i-1) -$

$\frac{(i+1)(i+2)}{2} + 1 + 2 = \frac{M^2}{2} - \frac{3M}{2} + 1 + i\left(2M - \frac{3}{2}\right) - \frac{i^2}{2} \quad (*)$

Then the median is that natural number i such that i is the first number for which $(*) > \frac{1}{2} \sum_{j=1}^{N-1} \alpha_j(K_N) = \frac{N(N-1)}{4} = M^2 - \frac{3}{2}M + \frac{1}{2}$. This means that

$$i^2 - (4M-3)i + (M^2 - 1) < 0$$

and thus $i = 2M - \frac{3}{2} \pm \sqrt{3M^2 - 6M + \frac{13}{4}}$. The plus sign is not possible as otherwise i > N, hence

$$i = Md_N > 2M - \frac{3}{2} - \sqrt{3M^2 - 6M + \frac{13}{4}}$$

$$> 2M - \frac{3}{2} - \sqrt{3} \, M = \left(2 - \sqrt{3}\right)M - \frac{3}{2}$$

Hence: $Md_{K_N} > \left(\frac{2-\sqrt{3}}{2}\right)N - \left(\frac{1+\sqrt{3}}{2}\right)$.

Consequently: $\lim_{N \to +\infty} \frac{Md_{K_N}}{\ln(N)} = +\infty$ which proves that $(K_N)_{N \in \mathbb{N}}$ is not(SWMd).

We further see that the median of the delta-sequence of $K_N$ (denoted here as $Md(K_N)$ ) is M-1 = (N-1)/2. Hence:

$$\lim_{N \to +\infty} \frac{M\delta_{K_N}}{\ln(N)} = \lim_{N \to +\infty} \frac{\frac{N-1}{2}}{\ln(N)} = +\infty$$

which shows that the sequence $(K_N)_{N \in \mathbb{N}}$ is DSWMd. □

Corollary. If $(\Omega_N)_{N \in \mathbb{N}}$ is DSWL then it is not necessarily SWD.

Proof. Assume that if $(\Omega_N)_{N \in \mathbb{N}}$ is DSWL then it is also SWD. Now DSWMd implies DSW (Proposition 1) from which we would know that $(\Omega_N)_{N \in \mathbb{N}}$ is SWD, from which it would follow by (Egghe & Rousseau, 2024) that $(\Omega_N)_{N \in \mathbb{N}}$ were SWMd, which is a contradiction (by Theorem 3).

Proposition 2. Consider $Z_1 = \{(\Omega_N)_{N \in \mathbb{N}} ; (\Omega_N)_{N \in \mathbb{N}} \; is \; SWD \}$ and $Z_2 = \{(\Omega_N)_{N \in \mathbb{N}} ; (\Omega_N)_{N \in \mathbb{N}} \; is \; DSWL \}$, then $Z_1 \cap Z_2 \neq \emptyset$.

Proof. It suffices to give one sequence $(\Omega_N)_{N \in \mathbb{N}}$ in the intersection $Z_1 \cap Z_2$. We know (Egghe & Rousseau, 2.6.2) that the sequence of stars is SWD and we also know that this sequence is DSWL. This proves this proposition.

Proposition 3. Consider $Z_3 = \{(\Omega_N)_{N \in \mathbb{N}} ; (\Omega_N)_{N \in \mathbb{N}} \; is \; not \; SWD \}$ and $Z_4 = \{(\Omega_N)_{N \in \mathbb{N}} ; (\Omega_N)_{N \in \mathbb{N}} \; is \; not \; DSWL \}$, then $Z_3 \cap Z_4 \neq \emptyset$.

Proof. Again it suffices to give one sequence in the intersection. The sequence $(CH_N)_{N \in \mathbb{N}}$ of N-chains is situated in the intersection, see Example 4.4.2. Note that also the sequence of N-polygons provides another example, see Example 4.4.3.

# 7. Delta-sequences and small worlds derived from degree distributions

## 7.1 Theorem 4

Consider the network sequences $(\Omega_N)_{N \in \mathbb{N}}$ and $(\Omega'_N)_{N \in \mathbb{N}}$ , such that for each $N \in \mathbb{N}$ , $\#\Omega_N = \#\Omega'_N$ . If now, there exists $N_0 \in \mathbb{N}$, such that for each $N \geq N_0 : \Delta_{\Omega_N}$ -< $\Delta_{\Omega'_N}$ , then

(a) $(\Omega_N)_{N \in \mathbb{N}}$ is DSWL implies that $(\Omega'_N)_{N \in \mathbb{N}}$ is DSWL, and hence DSWA.

(b) $(\Omega_N)_{N \in \mathbb{N}}$ is DSWA implies that $(\Omega'_N)_{N \in \mathbb{N}}$ is DSWA

(c) it does not follow that $(\Omega_N)_{N \in \mathbb{N}}$ is DSWMδ implies that $(\Omega'_N)_{N \in \mathbb{N}}$ is DSWMδ

(d) the opposite relations of (a) and (b) do not hold.

Proof. Results (a) and (b) follow from the definitions of DWDL and DSWA and Proposition 1 in 4.3. The opposite relations of (a) and (b) do not hold, because if they did then we would have an equivalence in that proposition, which does not hold.

Finally, we prove part (c). Inspired by the spider $S_M$ we construct the following networks. We consider three positive natural numbers M, a, and b (a and b stay fixed) and construct two networks with N = 2M+a+b nodes. For the first one, denoted as $S_{1,N}$, we take b < a. It consists of a complete (M+a) network, where, moreover, on (M+b) of these nodes we add one node (by a single link), see Fig. 12.

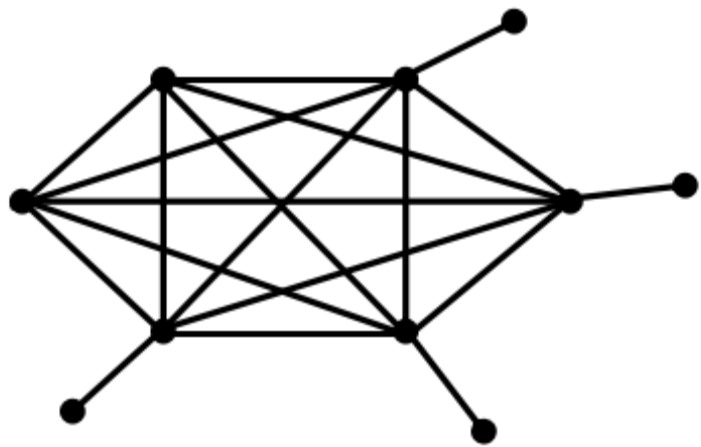

Fig. 12 Case M = 3, a = 3, b= 1

Then $\Delta_{S_{1,N}} = \left( \underbrace{M+a, \dots, M+a}_{(M+b)\,times}, \underbrace{M+a-1, \dots, M+a-1}_{(a-b)\,times}, \underbrace{1, \dots, 1}_{(M+b)\,times} \right)$.

For the second network, denoted as $S_{2,N}$, we take b > a. It again consists of a complete M+a network, on each of these nodes we add a singly-linked node, while moreover on b-a nodes we add a second, single-linked node, see Fig. 13.

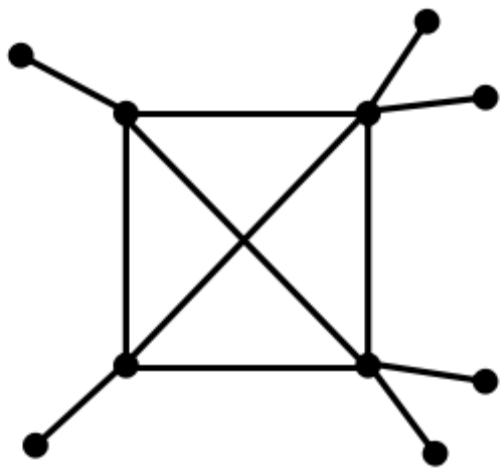

Fig. 13. Case M =3, a = 1, b=3

Then $\Delta_{S_{2,N}} = \left( \underbrace{M+a+1,\dots,M+a+1}_{(b-a)\,times}, \underbrace{M+a,\dots,M+a}_{(M+2a-b\,)\,times}, \underbrace{1,\dots,1}_{(M+b)\,times} \right).$

Clearly, $\Delta_{S_{1,N}} -\prec \Delta_{S_{2,N}}$ for each N. Now $M\delta(S_{1,N}) =$ M+a-1, which tends to infinity for N (or M) $\to \infty$. This shows that $\left(S_{1,N}\right)_N$ is DSWMd. Further, $M\delta\left(S_{2,N}\right) = 1$ (as b > a), which shows that $\left(S_{2,N}\right)_N$ is NOT(DSWMd).

## 7.2 Definition

Based on Theorem 4, we propose the following definition:

Given two networks G and H with delta-sequences $\Delta_G$ and $\Delta_H$, such that $\Delta_G -\prec \Delta_H$ then we say that network H is a smaller world than network G in the degree majorization sense.

We explain this statement: one cannot say that a network is a small world (at least not in our view as we defined the notion of a small world only for network sequences), but it is possible to compare two networks in terms of small-worldliness. This is done based on the generalized Lorenz curve. This is similar to the use of the classical Lorenz curve for comparing e.g., income inequality.

Theorem 4 shows that the uses of the terms "small world", "smaller world" and "small-worldliness" are consistent.

## 8. The neighboring array and the neighboring index

## 8.1 Definitions

Let G be a connected network with N nodes. We define the neighboring (or gamma) array of G, denoted as $\Gamma(G) = (\gamma_1(G), \gamma_2(G), \ldots, \gamma_N(G))$ , with $\gamma_i(G) = \sum_{j \in A(i)} \delta_j(G)$, with

$A(i) = \{j \in \{1, \ldots, N\}, i = j \text{ or there is a link between nodes i and j}\}$ .

Stated otherwise, the gamma-value of a node in a network is equal to the sum of the degree centralities of its zeroth and first-order neighbors. Next we define $v(G) = \sum_{i=1}^{N} \gamma_i(G)$ as the neighboring index of G.

## 8.2 A characterization of the neighboring array.

Let A be the adjacency matrix of an N-node network and let e = $\left(\underbrace{1, \ldots, 1}_{N \; times}\right)$ be the unit array. Then we have the following matrix multiplication result.

Theorem 5. $\Gamma = e.(A^2 + A)$

Proof. It is easy to check (and well-known) that $\Delta = e.A$. It is also well-known that the elements (i,j) of matrix $A^2$, denoted as $b_{ij}$ , yield the number of paths from node i to node j with length 2. Then $e.A^2 = \left(\sum_{j=1}^{N} b_{ij}\right)_{i=1,\ldots,N} = \left(\sum_{j=1}^{N} b_{ji}\right)_{i=1,\ldots,N} = \Gamma - \Delta$ □

## 8.3 Proposition 4

$$v(G) = \sum_{j=1}^{N} \delta_j(\delta_j + 1)) \qquad (14)$$

Proof. Every value $\delta_j$ occurs $\delta_j$ times for its first-order neighbors, plus one more time for its zeroth-order neighbor (itself). □

Remark. By (14) $\upsilon(G)$ follows from $\Delta(G)$, but when $\Gamma(G)$ is given, it is still an open problem if it is possible to construct $\Delta(G)$.

8.4 Examples (with N nodes); we assume that gamma-values are given in decreasing order.

8.4.1 The complete N-node network

$$\Gamma = \left( \underbrace{N(N-1), N(N-1), \dots, N(N-1)}_{N\ times} \right)$$

$$\upsilon = N^2(N-1)$$

8.4.2 The star

$$\Gamma = \left( 2(N-1), \underbrace{N, N, \dots, N}_{(N-1)\ times} \right)$$

$$\upsilon = (N+2)(N-1)$$

8.4.3 The polygon (N>2)

$$\Gamma = \left( \underbrace{6,6,\dots,6}_{N\ times} \right)$$

$$\upsilon = 6N$$

8.4.4 Chain (N>3)

$$\Gamma = \left( \underbrace{6,6,\dots,6}_{(N-4)times}\ 5,5,3,3 \right)$$

$$\upsilon = 6(N-4) + 16 = 6N - 8$$

8.4.5 The non-isomorphic networks (N=6).

Consider the networks shown in Fig. 14. Their alpha-sequences are the same, namely $(10, 5, 0, 0, 0)$ and so are their delta-sequences: $(4,4,3,3,3,3)$ , but their gamma-sequences are different, showing that these networks are not isomorphic: $\Gamma_G = (17,17,14,14,13,13)$ while $\Gamma_{G'} = (16,16,14,14,14,14)$ . Note that $\upsilon(G) = \sum_{i=1}^{N} \gamma_i = \sum_{i=1}^{N} \gamma_i' = \upsilon(G') = 88$.

This example shows that gamma-sequences are stronger than the combination of alpha- and delta-sequences when it comes to detecting isomorphisms. Yet, the combination of alpha-, delta- and gamma-sequences is not enough to detect isomorphism as shown in the next example.

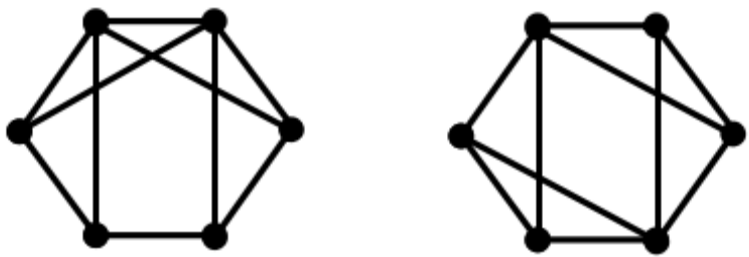

Fig. 14 Two non-isomorphic networks G (left) and G′ (right)

8.4.6 An example of two non-isomorphic networks with equal alpha-, delta-, and gamma-sequences.

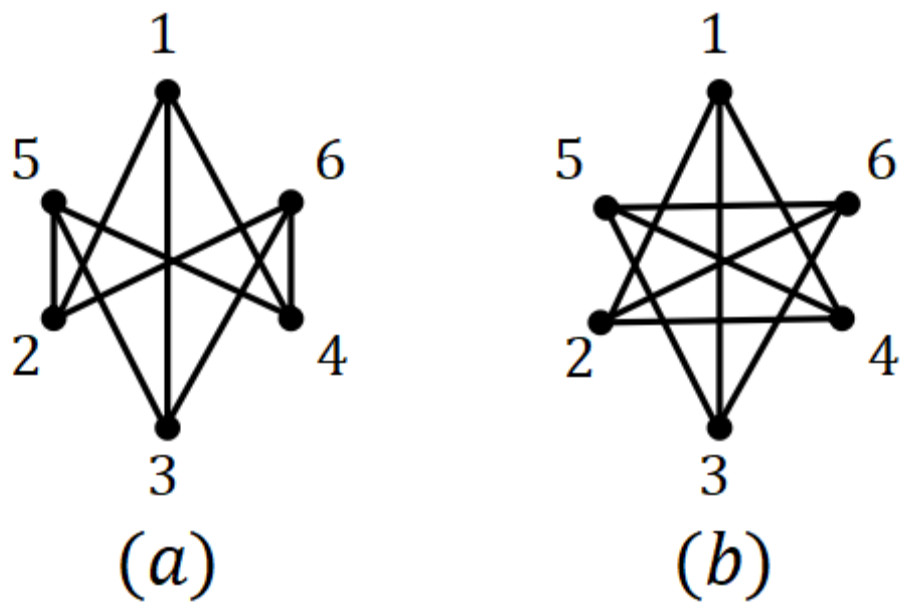

Fig.15. Two non-isomorphic networks with equal alpha-, delta- and gamma-sequences

For both networks in Fig. 15 we have AF=(9,6,0,0); $\Delta$ = (3,3,3,3,3,3) and $\Gamma$ = (12,12,12,12,12,12). Yet, they are not isomorphic as Fig. 15(a) has triples that are not connected, such as {1,5,6}, while such triples do not exist in Fig.15(b).

Remarks

(a) If two networks have the same total degree $\sum_{j=1}^{N} \delta_j$ then they do not necessarily have the same neighboring index $\sum_{j=1}^{N} \gamma_j$.

The following networks, shown in Fig.16 (note that they are trees) have the same $\sum_{j=1}^{N} \delta_j = 8$, but the chain has $\sum_{j=1}^{N} \gamma_j = 22$, while the other one has $\sum_{j=1}^{N} \gamma_j = 24$.

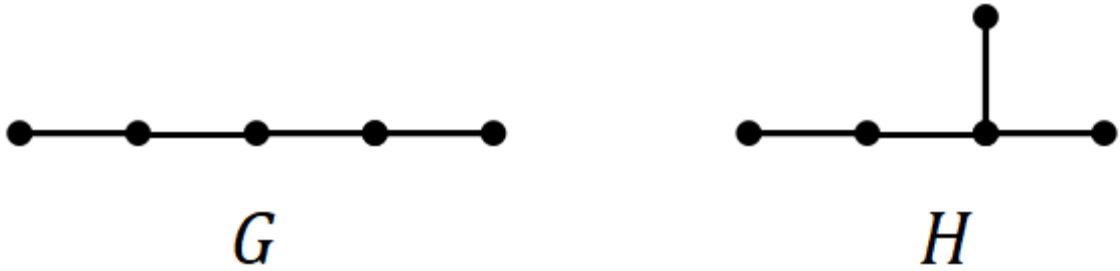

Fig. 16. Networks with the same total degree but different neighboring index (same example as in Fig. 1)

(b)   If two networks have the same neighboring index $\sum_{j=1}^{N} \gamma_j$ then they do not necessarily have the same total degree $\sum_{j=1}^{N} \delta_j$.

The networks shown in Fig. 17 (with N = 6) have the same neighboring index $\sum_{j=1}^{N} \gamma_j$, namely 48, but different total degrees $\sum_{j=1}^{N} \delta_j$, namely 14 and 12.

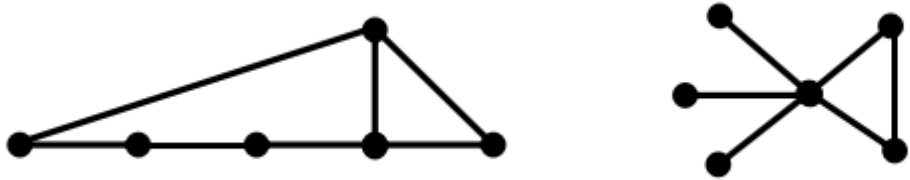

Fig. 17. Two networks with the same neighboring index but different total degree

We close this section by providing examples that the majorization relation -< is not kept between Δ and Γ.

Concretely:  Δ -< Δ' ⇏ Γ -< Γ' nor Γ' -< Γ                    (15)

and Γ -< Γ' ⇏ Δ -< Δ'  nor  Δ' -< Δ                    (16)

Indeed, for (15) we consider the networks (N=6), see Fig. 18.

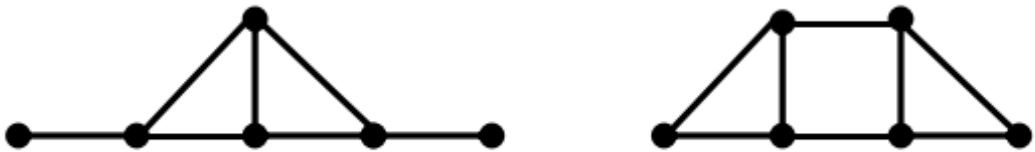

Fig.18. Illustration for (15)

The upper network has Δ = (3,3,3,3,1,1) and Γ = (12,12,10,10,4,4) while the lower one has Δ′ = (3,3,3,3,2,2) and Γ′ = (11,11,11,11,8,8). Then Δ -< Δ′, but neither Γ -< Γ′ nor Γ′-< Γ holds.

For the case (16) (with N =4) we consider Fig.19.

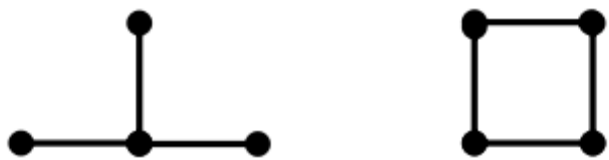

Fig. 19. Illustration for (16)

For the upper network, we have Δ = (3,1,1,1) and Γ = (6,4, 4,4) while the lower one has Δ′ = (2,2,2,2) and Γ′ = (6,6,6,6). Consequently: Γ -< Γ′, but neither Δ -< Δ′, nor Δ′ -< Δ  holds.

## 9. For trees, equal neighboring arrays lead to equal delta-arrays

### 9.1 Theorem 6

If G is a tree and $\Gamma(G) = \Gamma(G')$ then also G′ is a tree and $\Delta(G) = \Delta(G')$.

Proof. Assume that $\Delta(G) \neq \Delta(G')$. Then we know by the definition of the delta-sequence (1) that there exists a $\in \mathbb{N}$ such that

$$\#\left\{j \, ; \, \delta_j = a\right\} \neq \#\left\{j \, ; \, \delta_j' = a\right\}$$

Because of the inequality above we know that there exists $i \in \left\{j \, ; \, \delta_j = a\right\}$ and $i \notin \left\{j \, ; \, \delta_j' = a\right\}$ hence $\delta_i \neq \delta_i'$ and $\gamma_i = \gamma_i'$. As $\Gamma(G) = \Gamma(G')$ it follows that $\exists k \neq i$, node k directly linked to node i such that $\delta_k \neq \delta_k'$ and $\gamma_k = \gamma_k'$. Next, we continue by induction.

Assume we have different nodes $i_1, \ldots, i_m$ such that $\forall j = 1, \ldots, m$: $\delta_{i_j} \neq \delta_{i_j}'$ and $\gamma_{i_j} = \gamma_{i_j}'$, where, moreover, $\forall j = 2, \ldots, m$: node $i_j$ is directly linked to node $i_{j-1}$. As $\delta_{i_m} \neq \delta_{i_m}'$ and $\Gamma(G) = \Gamma(G')$ it follows that there exists a node $i_{m+1}$ directly linked to node $i_m$ such that $\delta_{i_{m+1}} \neq \delta_{i_{m+1}}'$ and $\gamma_{i_{m+1}} = \gamma_{i_{m+1}}'$. Now, node $i_{m+1}$ is not equal to any of the nodes $i_1, \ldots, i_m$ as otherwise there would exist a cycle in the tree G, which is impossible. Of course, networks G and G′ are are finite (like all networks in this article), which leads to a contradiction. Hence $\Delta(G) = \Delta(G')$ and by Knuth's Lemma, both networks are trees. □

## 10 Conclusion

We examined the degree distribution of a network and presented several examples. Utilizing the non-normalized Lorenz curve, we employed a generalized form of the majorization partial order. It is important to highlight that this represents a novel and fundamental application of the generalized Lorenz partial order. Our investigations support Rousseau's statement that Lorenz curves and Gini indices are universal tools for studying networks. Depending on the aim of the study the appropriate Lorenz-type curve should be used (Rousseau, 2011). Additionally, we introduced measures, including a Gini-type index, that respect the generalized Lorenz partial order.

We further introduced a new class of small worlds, namely those based on degree centralities of networks. Similar to a previous study, small worlds are defined as sequences of networks with certain limiting properties. We distinguish between three types of small worlds: those based on the highest degree, those based on the average degree, and those based on the median degree. We show that these new classes of small worlds are different from those introduced previously based on the diameter of the network or the average and median distance between nodes. However, there exist sequences of networks that qualify as small worlds in both senses of the word, with stars being an example. Our approach enables the comparison of two networks with an equal number of nodes in terms of their "small-worldliness". This comparison uses generalized Lorenz curves and the corresponding notion of generalized Lorenz majorization.

Extending the idea of delta- and alpha-sequences we introduced gamma-sequences, gave examples, showed their relation with delta-sequences, and showed that there exist non-isomorphic networks with the same alpha-, delta- and gamma-sequences.

We end this article by stating some open problems (OP):

OP1. Apply the generalized Lorenz order to the gamma-sequence.
OP2. Define and study Small Worlds in terms of the gamma-sequence.
OP3. Does the delta-sequence follow from the gamma-sequence?

Acknowledgment. The author thanks Li Li (Beijing) for drawing excellent illustrations, and Ronald Rousseau for stimulating discussions.